\documentclass[a4paper,reqno,10pt,english]{amsart}
\usepackage{geometry}

\usepackage{graphicx}
\graphicspath{{Figures/}}

\usepackage{times}
\usepackage{courier}
\usepackage[euler-digits]{eulervm}
\usepackage{eurosym}
\usepackage{hyperref}
\usepackage{amssymb,amsmath}

\usepackage{epsfig, graphicx}

\begin{document}

\title{Lotsize optimization leading to a $p$-median problem with cardinalities}

\author{Constantin Gaul}

\address{Constantin Gaul\\Department of Mathematics, Physic and Informatics\\University of Bayreuth\\Germany}
\email{cost.gaul@gmx.de}

\author{Sascha Kurz}

\address{Sascha Kurz\\Department of Mathematics, Physic and Informatics\\University of Bayreuth\\Germany}
\email{sascha.kurz@uni-bayreuth.de}

\author{J\"org Rambau}

\address{J\"org Rambau\\Department of Mathematics, Physic and Informatics\\University of Bayreuth\\Germany}
\email{j\"org.rambau@uni-bayreuth.de}

\begin{abstract}
We consider the problem of approximating the branch and size dependent demand of a fashion discounter with many branches
by a distributing process being based on the branch delivery restricted to integral multiples of lots from a small
set of available lot-types. We propose a formalized model which arises from a practical cooperation with an industry partner.
Besides an integer linear programming formulation and a primal heuristic for this problem we also consider a more abstract
version which we relate to several other classical optimization problems like the $p$-median problem, the facility location
problem or the matching problem.
\end{abstract}

\keywords{lotsize optimization, $p$-median problem, facility location problem, integer linear program formulation, primal heuristic, real world data, location-allocation}
\subjclass[2000]{90B80; 90C59, 90C10}

\maketitle

\section{Introduction}
\label{sec:Introduction}

Usually, fashion discounters can only achieve small profit margins.  Their economic success depends mostly in the ability to meet the customers' demands for individual products.  More specifically: offer exactly what you can sell to your customers.  This task has two aspects: offer what the customers would like to wear (attractive products) and offer the right volumes in the right places and the right sizes (demand consistent branch and size distribution).

In this paper we deal with the second aspect only: meet the branch and size specific demand for products as closely as possible. Our industry partner is a fashion discounter with more than 1\,000~branches most of whose products are never replenished, except for the very few {``}never-out-of-stock{''}-products (NOS products): because of lead times of around three months, apparel replenishments would be too late anyway.  In most cases the supplied items per product and apparel size lie in the range between $1$ and~$6$. Clearly there are some difficulties to determine a good estimate for the branch and size dependent demand, but besides a few practical comments on this problem we will blind out this aspect of the problem completely.

The problem we deal with in this article comes from another direction.  Our business partner is a discounter who has a lot of pressure to reduce its costs.  So he is forced to have a lean distribution logistics that works efficiently.  Due to this reason he, on the one hand, never replenishes and, on the other hand, tries to reduce the distribution complexity. To achieve this goal the supply of the branches is based on the delivery of lots, i.e., pre-packed assortments of single products in various sizes. Every branch can only be supplied with an integral multiple of one lot-type from a rather small number of available lot-types. So he has to face an approximation problem: which (integral) multiples of which (integral) lot-types should be supplied to a branch in order to meet a (fractional) mean demand as closely as possible?

We call this specific demand approximation problem the \emph{lot-type design problem (LDP)}.

\subsection{Related Work}
\label{subsec:RelatedWork}

The model we suggest for the LDP is closely related to the extensively
studied $p$-median- and the facility location problem. These problems
appear in various applications as some kind of clustering
problems. Loads of heuristics have been applied onto
them. Nevertheless the first constant-factor approximation algorithm,
based on LP rounding, was given not until 1999 by Charikar, Guha,
Tardos, and Shmoys \cite{ApproxAlgo2}.  We will give some more
detailed treatment or literature of approximation algorithms and
heuristics for the $p$-median- and the facility location problem in
Subsection \ref{subsec:ApproxAlgosAndHeuristics}.

\subsection{Our contribution}
\label{subsec:OurContribution}

In cooperation with our business partner, we identified the lot-type
design problem as a pressing real-world task.  We present an integer
linear program (ILP) formulation of the LDE that looks abstractly like
a $p$-median problem with an additional cardinality constraint.  We
call this problem the \emph{cardinality constrained $p$-median problem
  (Card-$p$-MP)}.  To the best of our knowledge, the Card-$p$-MP has
not been studied in the literature so far.

Although the ILP model can be solved by standard software on a
state-of-the-art PC in reasonable time, the computation times are
prohibitive for the use in the field, where interactive decision
support on a laptop is a must for negotiations with the supplier.
Therefore, we present a very fast primal any-time heuristics, that
yields good solutions almost instantly and searches for improvements
as long as it is kept running.  We demonstrate on real data that the
optimality gaps of our heuristics are mostly way below 1\,\%. At the
moment these heuristics are in test mode.

%% Starting from a practical problem we abstract the key properties and
%% formulate a new optimization problem, which we call cardinality
%% $p$-facility-location-problem. Relaxing this problem in different ways
%% gives some well known standard optimization problems like the
%% $p$-median-problem or the facility-location problem, which is also
%% known as the $p$-facility-location problem if those two problems are
%% combined.  The new enhancement of our problem is a cardinality
%% condition, which makes the problem somewhat harder. Besides an integer
%% linear program formulation we give a primal heuristic to solve the
%% problem very quickly on real world data. To demonstrate the potential
%% of our proposed primal heuristic we perform some test calculations to
%% determine the optimality gap on some realistic data sets.

\subsection{Outline of the paper}
\label{subsec:Outline}

In Section \ref{sec:RealWorldProblem} we will briefly describe the
real world problem, which we will formalize and model in Section
\ref{sec:Modeling}. In Section \ref{sec:CardPFacilityLocationProblem}
we will present its abstract version, the \emph{cardinality
  constrained $p$-median problem (Card-$p$-MP)}. Besides a formalized
description we relate it to several other well known optimization
problems like the matching problem, the facility location problem, or
the $p$-median problem. In Section \ref{sec:Heuristics} we present a
primal heuristic for the Card-$p$-MP, which we apply onto our real
world problem. We give some numerical data on the optimality gap of
our heuristic before we draw a conclusion in
Section~\ref{sec:Conclusion}.

\section{The real world problem}
\label{sec:RealWorldProblem}
Our industry partner is a fashion discounter with over $1\,000$
branches.  Products can not be replenished, and the number of sold
items per product and branch is rather small. There are no historic
sales data for a specific product available, since every product is
sold only for one sales period. The challenge for our industry partner
is to determine a suitable total amount of items of a specific product
which should be bought from the supplier. For this part the knowledge
and experience of the buyers employed by a fashion discounter is
used. We seriously doubt that a software package based on historic
sales data can do better.

But there is another task being more accessible for computer aided
forecasting methods. Once the total amount of sellable items of a
specific product is determined, one has to decide how to distribute
this total amount to a set of branches $B$ in certain apparel sizes
with in general different demands.  There are some standard techniques
how to estimate branch- and size-dependent demand from historic sales
data of related products, being, e.g., in the same commodity group.
We will address the problem of demand forecasting very briefly in
Subsection~\ref{subsec:Forecast}.  But let us assume for simplicity
that we know the exact (fractional) branch and size dependent mean
demands for a given new product or have at least good estimates.

Due to cost reasons, our industry partner organizes his distribution
process for the branches using a central warehouse. To reduce the
number of necessary handholds in the distributing process he utilizes
the concept of lots, by which we understand a collection of some items
of one product. One could have in mind different sizes or different
colors at this point. To reduce the complexity of the distribution
process also the number of used lot-types, e.g., different collections
of items, is limited to a rather small number. 

One could imagine that the branch- and size-dependent demand for a
specific product may vary broadly over the large set of branches. This
is at least the case for the branches of our industry partner. The
only flexibility to satisfy the demand in each single branch is to
choose a suitable lot-type from the small sets of available lot-types
and to choose a suitable multiplier, i.e., how many lots of a chosen
lot-type a specific branch should get. One should keep in mind that we
are talking about small multipliers here, i.e., small branches will
receive only one lot, medium sized branches will receive two lots, and
very big branches will receive three lots of a lot-type with, say, six
items.

The cost reductions by using this lot-based distribution system are
paid with a lack of possibility to approximate the branch and
size-dependent demand. So one question is, how many different lot-types
one should allow in order to be able to approximate the branch- and
size-dependent demand of the branches up to an acceptable deviation on
the one hand and to avoid a complex and cost intensive distribution
process in the central warehouse on the other hand. But also for a
fixed number of allowed lot-types the question of the best possible
approximation of the demand by using a lot-based supply of the
branches arises. In other words we are searching for an optimal
assignment of branches to lot-types together with corresponding
multipliers so that the deviation between the theoretical estimated
demand and the planned supply with lots is minimal. This is the main
question we will focus on in this paper.

\section{Mathematical modeling of the problem}
\label{sec:Modeling}

In this section we will prescind the real world problem from the
previous section and will develop an formulation as a well defined
optimization problem. Crucial and very basic objects for our
considerations are the set of branches $\mathcal{B}$, the set of sizes
$\mathcal{S}$ (in a more general context one could also think of a set
of variants of a product, like, e.g., different colors), and the set of
products~$\mathcal{P}$.

In practice, we may want to sell a given product $p\in\mathcal{P}$
only in some branches $\mathcal{B}_p\subseteq\mathcal{B}$ and only in
some sizes $\mathcal{S}_p\subseteq\mathcal{S}$ (clearly there are
different types of sizes for, e.g., skirts or socks). To model the
demand of a given branch $b\in\mathcal{B}_p$ for a given product
$p\in\mathcal{P}$ we use the symbol $\eta_{b,p}$, by which we
understand a mapping $\varphi_{b,p}$ from the set of sizes
$\mathcal{S}_p$ into a suitable mathematical object. This object may
be a random variable or simply a real number representing the mean
demand. In this paper we choose the latter possibility. For the sake
of a brief notation we regard $\eta_{b,p}$ as a vector
$\begin{pmatrix}\varphi_{b,p}\left(s_{i_1}\right)&\varphi_{b,p}\left(s_{i_2}\right)&\dots&\varphi_{b,p}\left(s_{i_r}\right)\end{pmatrix}\in\mathbb{R}^r$,
where we assume that $\mathcal{S}=\{s_1,\dots,s_t\}$ and
$\mathcal{S}_p=\{s_{i_1},\dots,s_{i_r}\}$ with $i_j<i_{j+1}$ for all
$j\in\{1,\dots,r-1\}$.

\subsection{Estimation of the branch- and size-dependent demand}
\label{subsec:Forecast}
For the purpose of this paper, we may assume that the demands
$\eta_{b,p}$ are given, but, since this is a very critical part in
practice, we would like to mention some methods how to obtain these
numbers. Marketing research might be a possible source. Another
possibility to estimate the demand for a product is to utilize
historic sales information. We may assume that for each product $p$
which was formerly sold by our retailer, each branch
$b\in\mathcal{B}$, each size $s\in\mathcal{S}$ and each day of sales
$d$ we know the number $\tau_{b,p}(d,s)$ of items which where sold in
branch $b$ of product $p$ in size $s$ during the first $d$ days of
sales. Additionally we assume, that we have a set
$\mathcal{U}\subseteq\mathcal{P}$ of formerly sold products which are
in some sense similar (one might think of the set of jeans if our new
product is also a jeans) to the new product $\tilde{p}$. By
$\mathcal{U}_{b,s}$ we denote the subset of products in $\mathcal{U}$,
which were traded by a positive amount in size $s$ in branch $b$ and
by $\chi_{b,s}(p)$ we denote a characteristic function which equals
$1$ if product $p$ is distributed in size $s$ to branch $b$, and
equals $0$ otherwise. For a given day of sales $d$ the value
\begin{equation}
  \label{eq_DemandEstimation}
  \tilde{\eta}_{b,\tilde{p},d}(s):=\frac{c}{|\mathcal{U}_{b,s}|}
  \sum_{u\in\mathcal{U}_{b,s}}\frac{\tau_{b,u}(d,s)\cdot\sum\limits_{b'\in\mathcal{B}_{\tilde{p}}}
  \sum\limits_{s'\in\mathcal{S}_{\tilde{p}}} \chi_{b',s'}(u)}{\sum\limits_{b'\in\mathcal{B}_{\tilde{p}}}
  \sum\limits_{s'\in\mathcal{S}_{\tilde{p}}}\tau_{b',u}(d,s')}
\end{equation}
might be a useable estimate for the demand $\eta_{b,\tilde{p}}(s)$,
after choosing a suitable scaling factor $c\in\mathbb{R}$ so that the
total estimate demand
$$
  \sum_{b\in\mathcal{B}_{\tilde{p}}}\sum_{s\in\mathcal{S}_{\tilde{p}}} \tilde{\eta}_{b,\tilde{p},d}(s)
$$
over all branches and sizes equals the total requirements. We would
like to remark that for small days of sale $d$ the quality of the
estimate $\tilde{\eta}_{b,\tilde{p},d}(s)$ suffers from the fact that
the stochastic noise of the consumer behavior is to dominating and for
large $d$ the quality of the estimate suffers from the fact of
stockout-substitution.

There are parametric approaches to this problem in the literature
(like Poisson-type sales processes).  In the data that was available
to us, we could not verify the main assumptions of such models,
though (not even close).

In our real world data set we have observed the fact that the sales
period of a product (say, the time by which 80~\% of the supply is
sold) varies a lot depending on the product. This effect is due to the
attractiveness of a given product (one might think of two T-shirts
which only differ in there color, where one color hits the vogue and
the other color does not). To compensate this effect we have chosen
the day of sales~$d$ in dependence of the product
$u\in\mathcal{U}_{b,s}$. More precisely, we have chosen $d_u$ so that
in the first $d_u$ days of sales a certain percentage of all items of
product $u$ where sold out over all branches and sizes.

Another possibility to estimate the demand is to perform the
estimation for the branch-dependent demand aggregated over all sizes
and the size-dependent demand for a given branch separately.

More sophisticated methods of demand estimation from historic sales
based on small data sets are, e.g., described in
\cite{paper_rio,paper_feldversuch}. Also research results from
forecasting NOS (never-out-of-stock) items, see, e.g.,
\cite{f_handbook,kok_fisher,sec_forecasting_handbook} for some
surveys, may be utilized. Also quite a lot of software-packages for
demand forecasting a available, see \cite{software_survey} for an
overview.

\subsection{Supply of the branches by lots}
To reduce handling costs in logistic and stockkeeping our business
partner orders his products from its external suppliers in so called
lots. These are assortments of several items of one product in
different sizes which form an entity. One could have a set of T-shirts
in different sizes in mind which are wrapped round by a plastic
foil. The usage of lots has the great advantage of reducing the number
of picks during the distribution process in a high-wage country like
Germany, where our partner operates.

Let us assume that the set of sizes for a given product $p$ is given
by $\mathcal{S}_p=\{s_{i_1},\dots,s_{i_r}\}$ with $i_j<i_{j+1}$ for
all $j\in\{1,\dots,r-1\}$. By a lot-type $l$ we understand a mapping
$\varphi:\mathcal{S}_p\rightarrow\mathbb{N}$, which can also be
denoted by a vector
$\begin{pmatrix}\varphi\left(s_{i_1}\right)&\varphi\left(s_{i_2}\right)&\dots&\varphi\left(s_{i_r}\right)\end{pmatrix}$
of non-negative integers.

By $\mathcal{L}$ we denote the set of applicatory lot-types. One could
imagine that a lot of a certain lot-type should not contain too many
items in order to be manageable. In the other direction it should also
not contain too few items in order to make use of the cost reduction
potential of the lot idea. Since the set of applicatory lot-types may
depend on a the characteristics of a certain product $p$ we specialize
this definition to a set $\mathcal{L}_p$ of manageable lot-types. (One
might imagine that a warehouseman can handle more T-shirts than, e.g.,
winter coats; another effect that can be modeled by a suitable set of
lot-types is to enforce that each size in $\mathcal{S}_p$ is supplied
to each branch in $\mathcal{B}_p$ by a positive amount due to
juridical requirements for advertised products.)

To reduce the complexity and the error-proneness of the distribution
process in a central warehouse, each branch $b\in\mathcal{B}_p$ is
supplied only with lots of one lot-type $l_{b,p}\in\mathcal{L}_p$. We
model the assignment of lot-types $l\in\mathcal{L}_p$ to branches
$b\in\mathcal{B}_p$ as a function
$\omega_p:\mathcal{B}_p\rightarrow\mathcal{L}_p$, $b\mapsto
l_{b,p}$. Clearly, this assignment $\omega_p$ is a decision variable
which can be used to optimize some target function. The only
flexibility that we have to approximate the branch-, size- and product
dependent demand $\eta_{b,p}$ by our delivery in lots is to supply an
integral multiple of $m_{b,p}$ items of lot-type $\omega_p(b)$ to branch
$b$. Again, we can denote this connection by a function
$m_p:\mathcal{B}_p\rightarrow\mathbb{N}$, $b\mapsto m_{b,p}$. Due to
practical reasons, also the total number
$\left|\omega_p\left(\mathcal{B}_p\right)\right|$ of used lot-types for
a given product is limited by a certain number $\kappa$.

\subsection{Deviation between supply and demand}
With the notation from the previous subsection, we can represent the
replant supply for branch $b$ with product $p$ as a vector
$m_p(b)\cdot \omega_p(b)\in\mathbb{N}^r$. To measure the deviation
between the supply $m_p(b)\cdot \omega_p(b)$ and the demand
$\eta_{b,p}$ we may utilize an arbitrary vector norm
$\Vert\cdot\Vert$. Mentionable vector norms in our context are the sum
of absolute values
$$
  \Vert\begin{pmatrix} v_1 & v_2 & \dots & v_r\end{pmatrix}\Vert_1:=\sum_{i=1}^r \left|v_i\right|,
$$
the maximum norm
$$
  \Vert\begin{pmatrix} v_1 & v_2 & \dots & v_r\end{pmatrix}\Vert_\infty:= \max\{|v_i|\,:\,1\le i\le r\},
$$
and the general $p$-norm
$$
  \Vert\begin{pmatrix} v_1 & v_2 & \dots & v_r\end{pmatrix}\Vert_p:= \sqrt[p]{\sum_{i=1}^r |v_i|^p}
$$
for real numbers $p>0$, which is also called the Euclidean norm for
$p=2$. With this we can define the deviation
$$
  \sigma_{b,l,m}:=\Vert \eta_{b,p}-m\cdot l\Vert_\star
$$
between demand $\eta_{b,p}$ and supply
$m\in\{1,\dots,M\}=:\mathcal{M}\subset\mathbb{N}$ times lot-type
$l\in\mathcal{L}_p$ for each branch $b\in\mathcal{B}_p$ and an
arbitrary norm $\Vert\cdot\Vert_\star$ for a given product
$p\in\mathcal{P}$. It depends on practical considerations which norm
to choose. The $\Vert\cdot\Vert_1$-norm is very insensitive in respect
to outliers in contrast to the $\Vert\cdot\Vert_\infty$-norm which is
absolutely sensitive with respect to outliers. A possible compromise
may be the Euclidean norm $\Vert\cdot\Vert_2$, but for most
considerations we choose the $\Vert\cdot\Vert_1$-norm because of its
robustness.  (We do not trust every single exact value in our demand
forecasts that much.)

For given functions $m_p$ and $\omega_p$ we can consider the deviation
vector
$$
  \Sigma_p:=\begin{pmatrix}\sigma_{b_1,\omega_p(b_1),m_p(b_1)} & \sigma_{b_2,\omega_p(b_2),m_p(b_2)} & \dots &
  \sigma_{b_q,\omega_p(b_q),m_p(b_q)}\end{pmatrix}
$$
if the set of branches is written as
$\mathcal{B}_p:=\{b_1,\dots,b_q\}$. To measure the total deviation of
supply and demand we can apply an arbitrary norm
$\Vert\cdot\Vert_\star$, which may be different from the norm to
measure the deviation of a branch, onto $\Sigma_p$. In this paper we
restrict ourselves on the $\Vert\cdot\Vert_1$-norm, so that we have
$$
  \Vert\Sigma_p\Vert_1=\sum_{b\in\mathcal{B}_p}\sigma_{b,\omega_p(b),m_p(b)}.
$$

\subsection{The cardinality condition}
For a given assignment $\omega_p$ of lot-types to branches and
corresponding multiplicities $m_p$ then quantity
$$
  I:=\sum_{b\in\mathcal{B}_p} m_p(b)\cdot\Vert\omega_p(b)\Vert_1\quad\in\mathbb{N}
$$
gives the total number of replant distributed items of product $p$
over all sizes and branches. From a practical point of view we
introduce the condition
\begin{equation}
  \label{eq_CardinalityCondition}
  \underline{I}\le I\le\overline{I},
\end{equation}
where $\underline{I}$, $\overline{I}$ are suitable integers. One might
imagine that our retailer may buy a part of already produced products
so that there is a natural upper bound $\overline{I}$ or that there
are some minimum quantities. Another interpretation may be that the
buying department of our retailer has a certain idea on the value of
$I$ but is only able to give an interval
$\left[\underline{I},\overline{I}\right]$. 

During our cooperation with our busines partner we have learned that
in practice you do not get what you order. If you order exactly $I$
items of a given product you will obtain $I$ plus minus some certain
percentage items in the end. (And their actually exists a certain
percentage up to which a retailer accepts a deviation between the
original order and the final delivery by its external suppliers as a
fulfilled contract.) 

Besides these and other practical reason to consider an interval
$\left[\underline{I},\overline{I}\right]$ for the total number of
items of a given product, there are very strong reasons not to replace
Inequalities~(\ref{eq_CardinalityCondition}) by an equation, as we
will explain in the following.  Let us consider the case where our
warehouse (or our external suppliers in a low-cost-country) is only
able to deal with a single lot-type per product. This is the case
$\kappa=1$. Let us further assume that there exists a rather small
integer $k$ (e.g. $k=20$) fulfilling $\Vert l\Vert_1\le k$ for all
$l\in\mathcal{L}_p$. If $I$ contains a prime divisor being larger than
$k$, then there exist no assignments multiplicities $m_p\in\mathbb{N}$
($\omega_p$ is a constant function due to $\kappa=1$) which lead to a
feasible solution of our problem. These number-theoretic influences
are somewhat ugly. In some cases the lead to the infeasibility of our
problem or to bad solutions with respect to the quality of the
demand-supply approximation in comparison to a relaxed version of the
problem, where the restrictions on $I$ are weaker. One could have in
mind the possibility of throwing one item into the garbage if this
will have a large impact on the quality of the demand-supply
approximation.

In Equation (\ref{eq_DemandEstimation}) for the demand estimation we
have used a certain number $\tilde{I}$ for the total number of items
to scale the demands $\eta_{b,p}$ by a factor $c$. From a more general
point of view it may also happen that the total demand
$$
  \sum_{b\in\mathcal{B}_p}\sum_{s\in\mathcal{S}_p} \eta_{b,p}(s)
$$
is not contained in the interval
$\left[\underline{I},\overline{I}\right]$. In this case the
$\Vert\cdot\Vert_1$-norm may not be very appropriate.  In our
estimation process, however, the demand forecasts in fact yield demand
percentages rather than absolute numbers.  The total volume is then
used to calculate the absolute (fractional) mean demand values, so
that in our work-flow the total demand is always in the target interval.

\subsection{The optimization problem}
\label{subsec:OptimizationProblem}
Summarizing the ideas and using the notations from the previous
subsections we can formulate our optimization problem in the following
form. We want to determine an assignment function
$\omega_p:\mathcal{B}_p\rightarrow\mathcal{L}_p$ and multiplicities
$m_p:\mathcal{B}_p\rightarrow\mathcal{M}=\{1,\dots,M\}\subset\mathbb{N}$
such that the total deviation between supply and demand
\begin{equation}
  \label{eq_total_deviation}
  \sum_{b\in\mathcal{B}_p}\sigma_{b,\omega_p(b),m_p(b)}
\end{equation}
is minimized with respect to the conditions
\begin{equation}
  \label{NumerOfLottypes}
  \left|\omega_p\left(\mathcal{B}_p\right)\right|\le\kappa
\end{equation}
and
\begin{equation}
  \label{CardinalityCondition}
  \underline{I}\,\,\le\,\,\sum_{b\in\mathcal{B}_p} m_p(b)\cdot\Vert\omega_p(b)\Vert_1\,\,\le\,\,\overline{I}.
\end{equation}
We use binary variables $x_{b,l,m}$, which are equal to $1$ if and
only if lot-type $l\in\mathcal{L}_p$ is delivered with multiplicity
$m\in\mathcal{M}$ to Branch $b$, and binary variables $y_l$, which are
$1$ if and only if at least one branch in $\mathcal{B}_p$ is supplied
with Lottype $l\in\mathcal{L}_p$.  With this, we can easily model out
problem as an integer linear program:
\begin{eqnarray}
  \label{OrderModel_Target}
  \min && \sum_{b\in\mathcal{B}_p}\sum_{l\in\mathcal{L}_p}\sum_{m\in\mathcal{M}} \sigma_{b,l,m}\cdot x_{b,l,m}\\
  \label{OrderModel_EveryBranchOneLottyp}
  s.t. && \sum_{l\in\mathcal{L}_p}\sum_{m\in\mathcal{M}} x_{b,l,m}=1\quad\quad\quad\forall b\in\mathcal{B}_p\\
  \label{OrderModel_Cardinality_1}
       && \sum_{b\in\mathcal{B}_p}\sum_{l\in\mathcal{L}_p}\sum_{m\in\mathcal{M}} m\cdot\Vert l\Vert_1\cdot x_{b,l,m}\le\overline{I}\\
  \label{OrderModel_Cardinality_2}
       && \sum_{b\in\mathcal{b}_p}\sum_{l\in\mathcal{L}_p}\sum_{m\in\mathcal{M}} m\cdot\Vert l\Vert_1\cdot x_{b,l,m}\ge\underline{I}\\
  \label{OrderModel_Binding}
       && \sum_{m\in\mathcal{M}} x_{b,l,m}\le y_l\quad\quad\quad \forall b\in\mathcal{B}_p\,\forall l\in\mathcal{L}_p\\
  \label{OrderModel_UsedLottypes}
       && \sum_{l\in\mathcal{L}_p} y_l\le \kappa\\
       && x_{b,l,m}\in\{0,1\}\quad\quad\quad\forall b\in\mathcal{B}_p\,\forall l\in\mathcal{L}_p\,\forall m\in\mathcal{M}\\
       && y_l\in\{0,1\}\quad\quad\quad\quad\,\,\,\,\forall l\in\mathcal{L}_p
\end{eqnarray}

The objective function (\ref{OrderModel_Target}) represents the sum
(\ref{eq_total_deviation}), since irrelevant tuples $(b,l,m)$ may be
downtroddened by $x_{b,l,m}=0$. Condition
(\ref{OrderModel_EveryBranchOneLottyp}) states that we assign for each
Branch~$b$ exactly one lot-type with a unique multiplicity. The
cardinality condition (\ref{CardinalityCondition}) is modeled by
Conditions (\ref{OrderModel_Cardinality_1}) and
(\ref{OrderModel_Cardinality_2}) and the restriction
(\ref{NumerOfLottypes}) on the number of used lot-types is modeled by
Condition (\ref{OrderModel_UsedLottypes}). The connection between the
$x_{b,l,m}$ and the $y_l$ is fixed in the usual Big-M condition
(\ref{OrderModel_Binding}). We would like to remark that the
LP-relaxation of this ILP formulation is very strong above all in
comparison to the more direct ILP formulation, where we assume the
branch deviation between supply and demand is measured by the
$\Vert\cdot\Vert_1$-norm:
\begin{eqnarray*}
  \min && \sum_{b\in\mathcal{B}_p}\sum_{s\in\mathcal{S}_p} z_{b,s}\\
  s.t. && \eta_{b,p}(s)-\alpha_{b,s} \le z_{b,s}\quad\quad\quad\forall b\in\mathcal{B}_p\,\forall s\in\mathcal{S}_p\\
       && \alpha_{b,s}-\eta_{b,p}(s) \le z_{b,s}\quad\quad\quad\forall b\in\mathcal{B}_p\,\forall s\in\mathcal{S}_p\\
       && \sum_{l\in\mathcal{L}_p}\sum_{m\in\mathcal{M}} x_{b,l,m}=1\quad\quad\quad\forall b\in\mathcal{B}_p\\
       && \sum_{b\in\mathcal{B}_p}\sum_{l\in\mathcal{L}_p}\sum_{m\in\mathcal{M}} m\cdot\Vert l\Vert_1\cdot x_{b,l,m}\le\overline{I}\\
       && \sum_{b\in\mathcal{B}_p}\sum_{l\in\mathcal{L}_p}\sum_{m\in\mathcal{M}} m\cdot\Vert l\Vert_1\cdot x_{b,l,m}\ge\underline{I}\\
       && \sum_{m\in\mathcal{M}} x_{b,l,m}\le y_l\quad\quad\quad \forall b\in\mathcal{B}_p\,\forall l\in\mathcal{L}_p\\
       && \sum_{l\in\mathcal{L}_p} y_l\le \kappa\\
       && \sum_{l\in\mathcal{L}_p}\sum_{m\in\mathcal{M}}m\cdot l[s]\cdot x_{b,l,m}=\alpha_{b,s}\quad\quad\quad
          \forall b\in\mathcal{B}_p\,\forall s\in\mathcal{S}_p\\
       && x_{b,l,m}\in\{0,1\}\quad\quad\quad\forall b\in\mathcal{B}_p\,\forall l\in\mathcal{L}_p\,\forall m\in\mathcal{M}\\
       && y_l\in\{0,1\}\quad\quad\quad\quad\,\,\,\,\forall l\in\mathcal{L}_p\\
       && \alpha_{b,s}\in\mathbb{R}_0^+\quad\quad\quad\quad\quad\forall b\in\mathcal{B}_p\,\forall s\in\mathcal{S}_p,
\end{eqnarray*}
where $l[s]$ is the entry in Vector~$l$ corresponding to Size~$s$.

We would like to remark that our strong ILP formulation of the problem
of Subsection \ref{subsec:OptimizationProblem} can be used to solve
all real world instances of our business partner in at most 30~minutes
by using a standard ILP solver like \texttt{CPLEX 11}. Unfortunately,
this is not fast enough for our real world application. The buyers of
our retailer need a software tool which can produce a near optimal
order recommendation in real time on a standard laptop. The buying
staff travels to one of the external suppliers to negotiate several
orderings. When they get to the details, the buyer inserts some key
data like $\underline{I}$, $\overline{I}$, $\mathcal{B}_p$,
$\mathcal{S}_p$, and $\mathcal{L}_p$ into his laptop and immediately
wants a recommendation for an order in terms of multiples of
lot-types. For this reason, we consider in Section \ref{sec:Heuristics}
a fast heuristic, which has only a small gap compared to the optimal
solution on a test set of real world data of our business partner.

\section{The Cardinality Constrained $p$-Median Problem}
\label{sec:CardPFacilityLocationProblem}

In the previous section we have modeled our real world problem from
Section \ref{sec:RealWorldProblem}. Now we want to abstract from this
practical problem and formulate a more general optimization problem
which we will relate to several well known optimization problems.

For the general Cardinality Constrained $p$-Median Problem let $p$ be
an integer, $\mathcal{S}$ a set of chooseable items, $\mathcal{D}$ a
set of demanders, a demand function
$\delta:\mathcal{D}\rightarrow\mathbb{R}^+$, and
$\left[\underline{I},\overline{I}\right]\subseteq\mathcal{N}$ an
interval. We are looking for an assignment
$\omega:\mathcal{D}\rightarrow\mathcal{S}$ with corresponding
multipliers $m:\mathcal{D}\rightarrow\mathbb{N}$, such that the sum of
distances
$$
  \sum_{d\in\mathcal{D}} \Vert \delta(d)-m(d)\cdot\omega(d)\Vert
$$
is minimized under the conditions
$$
  |\omega(\mathcal{D})|\le p
$$
and
$$
  \underline{I}\quad\le\quad\sum_{d\in\mathcal{D}} m(d)\cdot |\omega(d)|\quad\le\quad\overline{I}.
$$

Let us now bring this new optimization problem in line with known
combinatorial optimizations problems.  Since we have to choose an
optimal subset of $\mathcal{S}$ to minimize a cost function subject to
some constraints the cardinality constrained $p$-median problem belongs to the
large class of generic selection problems.

Clearly, it is closely related to the $p$-median problem. The only
characteristics of our problem that are not covered by the $p$-median
problem are the multipliers $m$ and the cardinality condition. If we
relax the cardinality condition we can easily transform our problem
into a classical $p$-median problem. For every element
$d\in\mathcal{D}$ and every element $s\in\mathcal{S}$ there exists an
optimal multiplier $m_{d,s}$ such that $\Vert \delta(d)-m_{d,s}\cdot
s\Vert$ is minimal.

If we do not bound $|\omega(\mathcal{D})|$ from above but assign costs
for using elements of $\mathcal{S}$ instead, which means using another
lot-type in our practical application, we end up with the
facility location problem. Clearly we also have some kind of an
assignment-problem, since the have to determine an assignment $\omega$
between the sets $\mathcal{D}$ and a subset of $\mathcal{S}$.

One can also look at our problem from a completely different
angle. Actually we are given a set of $|\mathcal{B}|$ real-valued
demand-vectors, which we want to approximate by a finite number of
integer-valued vectors using integral multiples. There is a well
established theory in number theory on so called Diophantine
approximation \cite{Diophant_approx2,Diophant_approx} or simultaneous
approximation, which is somewhat related to our approximation
problem. Here one is interested in simultaneously minimizing
$$
  \left\Vert\alpha_i-\frac{p_i}{q}\right\Vert
$$
for linearly independent real numbers $\alpha_i$ by integers $p_i$ and
$q$ \cite{simult_approx,simult_approx2}. One might use some results
from this theory to derive some bounds for our problem. One might also
have a look at \cite{sim_approx_combinatorial}.

For a more exhaustive and detailed analysis of the taxonomy of the
broad field of facility-location problems and their modeling we refer
to \cite{ApproxAlgo19}.

\subsection{Approximation algorithms and heuristics for related problems}
\label{subsec:ApproxAlgosAndHeuristics}

Facility location problems and the $p$-median problem are well known
and much research has been done. Since, moreover, these problems are
closely related to our optimization problem, we would like to mention
some literature and methods on approximation algorithms and heuristics
for these problems.

Lin and Vitter \cite{ApproxAlgo9} have developed a filtering and
rounding technique which rounds fractional solutions of the standard
LP for these problems to obtain good integer solution. For the metric
case some some bounds for approximation quality are given. Based on
this work some improvements were done in \cite{ApproxAlgo10}, were the
authors give a polynomial-time $3.16$-approximation algorithm for the
metric facility location problem, and \cite{ApproxAlgo3,ApproxAlgo2},
where the authors give a polynomial-time $\frac{20}{3}$-approximation
algorithm for the metric $p$-median problem and a $9.8$-approximation
algorithm for the $p$-facility location problem.

Besides Rounding techniques of LP-solutions also greedy techniques
have been applied to the facility location problem and the $p$-median
problems. Some results are given in
\cite{ApproxAlgo5,ApproxAlgo7,ApproxAlgo8}. Since these problems are
so prominent in applications the whole broadness of heuristics are
applied onto it. Examples are scatter search
\cite{ApproxAlgo4,ApproxAlgo16}, local search
\cite{ApproxAlgo1,ApproxAlgo13}, and neighborhood search
\cite{ApproxAlgo14,ApproxAlgo15}.

Good overviews for the broad topic of approximation algorithms and
heuristics for the facility location and the $p$-median problem are
given in \cite{ApproxAlgo10,ApproxAlgo11,ApproxAlgo17,ApproxAlgo18}.

Besides results for the metric case there are also results for the
non-metric case, see, e.g., \cite{ApproxAlgo12}.

Unfortunately, none of the theoretical guarantees seems to survive the
introduction  of the cardinality constraint in general.

\section{A practical heuristic for the Cardinality Constrained $p$-Median Problem}
\label{sec:Heuristics}

As already mentioned in Section \ref{sec:Modeling} solving our ILP
formulation of our problem is too slow in practical applications. So
there is a real need for a fast heuristic which yields good solutions,
which is the top of this section.

In Section \ref{sec:CardPFacilityLocationProblem} we have analyzed our
problem from different theoretical point of views. What happens if we
relax some conditions or fix some decisions. A very important decision
is: which lot-types should be used in the first place?  Here one should
have in mind that the cardinality $\left|\mathcal{L}_p\right|$ of the
set of feasible lot-types is very large compared to the number $\kappa$
of lot-types which can be used for the delivery process of a specific
product $p$.

\subsection{Heuristic selection of lot-types}
\label{subsec:SelectionLottypes}
For this selection problem of lot-types we utilize a scoring
method. For every branch $b\in\mathcal{B}_p$ with demand $\eta_{b,p}$
there exists a lot-type $l\in\mathcal{L}_p$ and a multiplicity
$m\in\mathbb{N}$ such that $\Vert\eta_{b,p}-m\cdot l\Vert$ is minimal
in the set $\left\{\Vert\eta_{b,p}-m'\cdot
  l'\Vert\,:\,l'\in\mathcal{L}_p,\,m'\in\mathbb{N}\right\}$. So for
every branch $b\in\mathcal{B}_p$ there exists a lot-type that fits
best. More general, for a given $k\le\left|\mathcal{L}_p\right|$ there
exist lot-types $l_1,\dots,l_k$ such that $l_i$ fits $i$-best if one
uses the corresponding optimal multiplicity. Let us examine this
situation from the point of view of the different lot-types. A given
lot-type $l\in\mathcal{L}_p$ is the $i$-best fitting lot-type for a number
$\varrho_{l,i}$ of branches in $\mathcal{B}_p$. Writting these numbers
$\varrho_{l,i}$ as a vector $\varrho_l\in\mathbb{N}^k$ we obtain score
vectors for all lot-types $l\in\mathcal{L}_p$.

Now we want to use these score vectors $\varrho_l$ to sort the
lot-types of $\mathcal{L}_p$ in decreasing \textit{approximation
  quality}. Using the lexicographic ordering $\preceq$ on vectors we
can determine a bijective rank function
$\lambda:\mathcal{L}_p\rightarrow\left\{1,\dots,\left|\mathcal{L}_p\right|\right\}$. (We
simply sort the score vectors according to $\preceq$ and for the case
of equality we choose an arbitrary succession.) We extend $\lambda$ to
subsets $\mathcal{L}'\subseteq\mathcal{L}_p$ by
$\lambda(\mathcal{L}')=\sum\limits_{l\in\mathcal{L}'}\lambda(l)\in\mathbb{N}$.

To fix the lot-types we simply loop over subsets
$\mathcal{L}'\subseteq\mathcal{L}_p$ of cardinality $\kappa$ in
decreasing order with respect to $\lambda(\mathcal{L}')$. In principle
we consider all possible selections $\mathcal{L}'$ of $\kappa$
lot-types, but in practise we stop our computations after a adequate
time period with the great advantage that we have checked the in some
heuristic sense most promising selections $\mathcal{L}'$ first.

Now we have to go into detail how to efficiently determine the $p$
best fitting lot-types with corresponding optimal multiplicities for
each branch $b\in\mathcal{B}_p$. We simply loop over all branches
$b\in\mathcal{B}_p$ and determine the set of the $p$ best fitting
lot-types separately. Here we also simply loop over all lot-types
$l\in\mathcal{L}_p$ and determine the corresponding optimal multiplier
$m$ by binary search (it is actually very easy to effectively
determine lower and upper bounds for $m$ from $\eta_{b,p}$ and $l$)
due to the convexity of norm functions. Using a heap data structure
the sorting of the $p$ best fitting lot-types can be done in
$O(|\mathcal{L}_p|)$ time if $k\log k\in O(|\mathcal{L}_p|)$, which is
not a real restriction for practical problems. We further want to
remark that we do not have to sort the score vectors completely since
in practice we will not loop over all
${\left|\mathcal{L}_p\right|}\choose{\kappa}$ possible selections of
lot-types. If one does not want to use a priori bounds (meaning that
one excludes the lot-types with high rank $\lambda$) one could use a
\textit{lazy} or delayed computation of the sorting of $\lambda$ by
utilizing again a heap data structure.

\subsection{Adjusting a delivery plan to the cardinality condition}
\label{subsec:AdjustingCardinality}
If we determine assignments $\omega_p$ with corresponding multipliers
$m_p$ with the heuristic being described in Subsection
\ref{subsec:SelectionLottypes} in many cases we will not satisfy the
cardinality condition (\ref{eq_CardinalityCondition}) since it is
totally unaccounted by our heuristic. Our strategy to satisfy the
cardinality condition (\ref{eq_CardinalityCondition}) is to adjust
$m_p$ afterwards by decreasing or increasing the calculated
multipliers unless condition (\ref{eq_CardinalityCondition}) is
fulfilled by pure chance.

Here we want to use a greed algorithm and have to distinguish two
cases. If $I(\omega_p,m_p)$ is smaller then $\underline{I}$, then we
increase some of the values of $m_p$, other wise we have
$I(\omega_p,m_p)>\overline{I}$ and we decrease some of the values of
$m_p$. Our procedure works iteratively and we assume that the current
multipliers are given by $\widetilde{m}_p$. Our stopping criteria is
given by $\underline{I}\le I(\omega_p,\widetilde{m}_p)\le\overline{I}$
or that there are no feasible operations left. We restrict our
explanation of a step of the iteration to the case where we want to
decrease the values of $\widetilde{m}_p$. For every branch
$b\in\mathcal{B}_p$ the reduction of $\widetilde{m}_p(b)$ by one
produces costs
$$
  \Delta_b^-=\sigma_{b,\omega_p(b),\widetilde{m}_p(b)-1}-\sigma_{b,\omega_p(b),\widetilde{m}_p(b)}
$$
if the reduction of $\widetilde{m}_p(b)$ by one is allowed (a suitable
condition is $\widetilde{m}_p\ge 1$ or $\widetilde{m}_p\ge 2$) and
$\Delta_b^-=\infty$ if we do not have the possibility to reduce the
multiplier $\widetilde{m}_p(b)$ by one. A suitable data structure for
the $\Delta_b^-$ values is a heap, for which the update after an
iteration can be done in $O(1)$ time. If we reach
$I(\omega_p,\widetilde{m}_p)<\underline{I}$ at some point, we simply
discard this particular selection $\omega_p$ and consider the next
selection candidate.

Since this adjustment step can be performed very fast one might also
take some kind of general swap techniques into account. Since for
these techniques there exists an overboarding amount of papers in the
literature we will not go into detail here, but we would like to
remark that in those cases (see Subsection \ref{subsec:OptimalityGap})
where the optimality gap of our heuristic lies above 1~\% swapping can
improve the solutions of our heuristic by a large part.

\subsection{Optimality gap}
\label{subsec:OptimalityGap}
To substantiate the usefullness of our heuristic we have compared the
quality of the solutions given by this heuristic after one second of
computation time (on a standard laptop) with respect to the solution
given by \texttt{CPLEX 11}.

Our business partner has provided us historic sales information for
nine different commodity groups each ranging over a sales period of at
least one and a half year. For each commodity group we have performed
a test calculation for $\kappa\in\{1,2,3,4,5\}$ distributing some
amount of items to almost all branches.

\bigskip

\noindent
\textbf{Commodity group 1}: % 20-4-4

\noindent
Cardinality interval: $[10630,11749]$\\
number of sizes: $|\mathcal{S}_p|=5$\\
number of branches: $|\mathcal{B}_p|=1119$

\begin{table}[htp]
  \begin{center}
    \begin{tabular}{|c|c|c|c|c|c|}
      \hline
      & $\kappa=1$ & $\kappa=2$ & $\kappa=3$ & $\kappa=4$ & $\kappa=5$ \\
      \hline
      \texttt{CPLEX} & 4033.34 & 3304.10 & 3039.28 & 2951.62 & 2891.96 \\
      \hline
      heuristic      & 4033.85 & 3373.95 & 3076.55 & 3011.49 & 2949.31 \\
      \hline
      gap            & 0.013\% & 2.114\% & 1.226\% & 2.028\% & 1.983\% \\
      \hline
    \end{tabular}
    \caption{Optimality gap in the $\Vert\cdot\Vert_1$-norm for our heuristic on commodity group 1}
    \label{table:CommodityGroup1}
  \end{center}
\end{table}

\noindent
\textbf{Commodity group 2}: % 2-3-3

\noindent
Cardinality interval: $[10000,12000]$\\
number of sizes: $|\mathcal{S}_p|=5$\\
number of branches: $|\mathcal{B}_p|=1091$

\begin{table}[htp]
  \begin{center}
    \begin{tabular}{|c|c|c|c|c|c|}
      \hline
      & $\kappa=1$ & $\kappa=2$ & $\kappa=3$ & $\kappa=4$ & $\kappa=5$ \\
      \hline
      \texttt{CPLEX} & 2985.48 & 2670.04 & 2482.23 & 2362.75 & 2259.57 \\
      \hline
      heuristic      & 3371.64 & 2671.72 & 2483.52 & 2362.90 & 2276.32 \\
      \hline
      gap            & 12.934\% & 0.063\% & 0.052\% & 0.006\% & 0.741\% \\
      \hline
    \end{tabular}
    \caption{Optimality gap in the $\Vert\cdot\Vert_1$-norm for our heuristic on commodity group 2}
    \label{table:CommodityGroup2}
  \end{center}
\end{table}

\noindent
\textbf{Commodity group 3}: % 2-7-4

\noindent
Cardinality interval: $[9785,10815]$\\
number of sizes: $|\mathcal{S}_p|=5$\\
number of branches: $|\mathcal{B}_p|=1030$

\begin{table}[htp]
  \begin{center}
    \begin{tabular}{|c|c|c|c|c|c|}
      \hline
      & $\kappa=1$ & $\kappa=2$ & $\kappa=3$ & $\kappa=4$ & $\kappa=5$ \\
      \hline
      \texttt{CPLEX} & 3570.3282 & 3022.2655 & 2622.8209 & 2488.1009 & 2413.55 \\
      \hline
      heuristic      & 3571.61 & 3023.91 & 2625.29 & 2492.07 & 2417.65 \\
      \hline
      gap            & 0.036\% & 0.054\% & 0.094\% & 0.160\% & 0.170\% \\
      \hline
    \end{tabular}
    \caption{Optimality gap in the $\Vert\cdot\Vert_1$-norm for our heuristic on commodity group 3}
    \label{table:CommodityGroup3}
  \end{center}
\end{table}

\noindent
\textbf{Commodity group 4}: % 3-3-2

\noindent
Cardinality interval: $[10573,11686]$\\
number of sizes: $|\mathcal{S}_p|=5$\\
number of branches: $|\mathcal{B}_p|=1119$

\begin{table}[htp]
  \begin{center}
    \begin{tabular}{|c|c|c|c|c|c|}
      \hline
      & $\kappa=1$ & $\kappa=2$ & $\kappa=3$ & $\kappa=4$ & $\kappa=5$ \\
      \hline
      \texttt{CPLEX} & 4776.36 & 4364.63 & 4169.94 & 4023.60 & 3890.87 \\
      \hline
      heuristic      & 5478.19 & 4365.47 & 4170.23 & 4024.55 & 3892.35 \\
      \hline
      gap            & 14.694\% & 0.019\% & 0.007\% & 0.024\% & 0.038\% \\
      \hline
    \end{tabular}
    \caption{Optimality gap in the $\Vert\cdot\Vert_1$-norm for our heuristic on commodity group 4}
    \label{table:CommodityGroup4}
  \end{center}
\end{table}

\noindent
\textbf{Commodity group 5}: % 6-3

\noindent
Cardinality interval: $[16744,18506]$\\
number of sizes: $|\mathcal{S}_p|=5$\\
number of branches: $|\mathcal{B}_p|=1175$

\begin{table}[htp]
  \begin{center}
    \begin{tabular}{|c|c|c|c|c|c|}
      \hline
      & $\kappa=1$ & $\kappa=2$ & $\kappa=3$ & $\kappa=4$ & $\kappa=5$ \\
      \hline
      \texttt{CPLEX} & 4178.71 & 3418.37 & 3067.74 & 2874.70 & 2786.69 \\
      \hline
      heuristic      & 4179.23 & 3418.87 & 3068.25 & 2875.21 & 2787.21 \\
      \hline
      gap            & 0.013\% & 0.015\% & 0.017\% & 0.018\% & 0.019\% \\
      \hline
    \end{tabular}
    \caption{Optimality gap in the $\Vert\cdot\Vert_1$-norm for our heuristic on commodity group 5}
    \label{table:CommodityGroup5}
  \end{center}
\end{table}

\noindent
\textbf{Commodity group 6}: % 21-5

\noindent
Cardinality interval: $[11000,13000]$\\
number of sizes: $|\mathcal{S}_p|=4$\\
number of branches: $|\mathcal{B}_p|=1030$

\begin{table}[htp]
  \begin{center}
    \begin{tabular}{|c|c|c|c|c|c|}
      \hline
      & $\kappa=1$ & $\kappa=2$ & $\kappa=3$ & $\kappa=4$ & $\kappa=5$ \\
      \hline
      \texttt{CPLEX} & 2812,22 & 2311,45 & 2100,78 & 1987,46 & 1909,21 \\
      \hline
      heuristic      & 2812,63 & 2311,87 & 2101,25 & 1987,93 & 1909,63 \\
      \hline
      gap            & 0.015\% & 0.018\% & 0.022\% & 0.024\% & 0.022\% \\
      \hline
    \end{tabular}
    \caption{Optimality gap in the $\Vert\cdot\Vert_1$-norm for our heuristic on commodity group 6}
    \label{table:CommodityGroup6}
  \end{center}
\end{table}

\noindent
\textbf{Commodity group 7}: % 20-5

\noindent
Cardinality interval: $[15646,17293]$\\
number of sizes: $|\mathcal{S}_p|=5$\\
number of branches: $|\mathcal{B}_p|=1098$

\begin{table}[htp]
  \begin{center}
    \begin{tabular}{|c|c|c|c|c|c|}
      \hline
      & $\kappa=1$ & $\kappa=2$ & $\kappa=3$ & $\kappa=4$ & $\kappa=5$ \\
      \hline
      \texttt{CPLEX} & 4501.84 & 3917.96 & 3755.20 & 3660.32 & 3575.55 \\
      \hline
      heuristic      & 4719.06 & 3918.46 & 3755.70 & 3660.84 & 3576.04 \\
      \hline
      gap            & 4.825\% & 0.013\% & 0.013\% & 0.014\% & 0.014\% \\
      \hline
    \end{tabular}
    \caption{Optimality gap in the $\Vert\cdot\Vert_1$-norm for our heuristic on commodity group 7}
    \label{table:CommodityGroup7}
  \end{center}
\end{table}

\noindent
\textbf{Commodity group 8}: % 2-1-2

\noindent
Cardinality interval: $[11274,12461]$\\
number of sizes: $|\mathcal{S}_p|=5$\\
number of branches: $|\mathcal{B}_p|=989$

\begin{table}[htp]
  \begin{center}
    \begin{tabular}{|c|c|c|c|c|c|}
      \hline
      & $\kappa=1$ & $\kappa=2$ & $\kappa=3$ & $\kappa=4$ & $\kappa=5$ \\
      \hline
      \texttt{CPLEX} & 3191.66 & 2771.89 & 2575.37 & 2424.31 & 2331.67 \\
      \hline
      heuristic      & 3579.35 & 2772.33 & 2575.81 & 2424.75 & 2332.11 \\
      \hline
      gap            & 12.147\% & 0.016\% & 0.017\% & 0.018\% & 0.019\% \\
      \hline
    \end{tabular}
    \caption{Optimality gap in the $\Vert\cdot\Vert_1$-norm for our heuristic on commodity group 8}
    \label{table:CommodityGroup8}
  \end{center}
\end{table}

\noindent
\textbf{Commodity group 9}: % 2-5-3

\noindent
Cardinality interval: $[9211,10181]$\\
number of sizes: $|\mathcal{S}_p|=5$\\
number of branches: $|\mathcal{B}_p|=808$

\begin{table}[htp]
  \begin{center}
    \begin{tabular}{|c|c|c|c|c|c|}
      \hline
      & $\kappa=1$ & $\kappa=2$ & $\kappa=3$ & $\kappa=4$ & $\kappa=5$ \\
      \hline
      \texttt{CPLEX} & 3616.71 & 3215.17 & 2981.02 & 2837.66 & 2732.29 \\
      \hline
      heuristic      & 3617.09 & 3215.53 & 3009.01 & 2860.85 & 2758.39 \\
      \hline
      gap            & 0.010\% & 0.011\% & 0.939\% & 0.817\% & 0.955\% \\
      \hline
    \end{tabular}
    \caption{Optimality gap in the $\Vert\cdot\Vert_1$-norm for our heuristic on commodity group 9}
    \label{table:CommodityGroup9}
  \end{center}
\end{table}

\bigskip

Besides these nine test calculations we have done several calculations
on our data sets with different parameters, we have, e.g., considered
case with fewer sizes, fewer branches, smaller or larger cardinality
intervals, larger $\kappa$, or other magnitudes for the cardinality
interval. The results are from a qualitative point of view more or
less the same, as for the presented test calculations.

\section{Conclusion and Outlook}
\label{sec:Conclusion}

Starting from a real world optimization problem we have formalized a
new general optimization problem, which we call cardinality
$p$-facility location problem. It turns out that this problem is
related to several other well known standard optimization problems. In
Subsection \ref{subsec:OptimizationProblem} we have given an integer
linear programming formulation which has a very strong
LP-relaxation. Nevertheless this approach is quit fast (computing
times below one hour), there was a practical need for fast heuristics
to solve the problem. We have presented one such heuristic which
performs very well on real world data sets with respect to the
optimality gap.

Some more theoretic work on the cardinality $p$-facility location
problem and its relationships to other classical optimization methods
may lead to even stronger integer linear programming formulations or
faster branch-and-bound frameworks enhanced with some graph theoretic
algorithms.

We leave also the question of a good approximation algorithm for the
cardinality $p$-facility location problem. Having the known
approximation algorithms for the other strongly related classical
optimization problems in mind, we are almost sure that it should be
not too difficult to develop good approximation algorithms for our
problem.

For the practical problem the uncertainties and difficulties
concerning the demand estimation have to be faced. There are several
ways to make solutions of optimization problems more robust. One
possibility is to utilize robust optimization methods.  Another
possibility is to consider the branch- and size dependent demands as
stochastic variables and to utilize integer linear stochastic
programming techniques. See, e.g., \cite{StochasticProgramming} or more
specifically \cite{ApproxAlgo11}.  These enhanced models, however,
will challenge the solution methods a lot, since the resulting
problems are of a much larger scale than the one presented in this
paper.  Nevertheless, this is exactly what we are looking at next.

\nocite{ApproxAlgo1,ApproxAlgo2,ApproxAlgo3,ApproxAlgo4,ApproxAlgo5,ApproxAlgo6,ApproxAlgo7,ApproxAlgo8,ApproxAlgo9,ApproxAlgo10,ApproxAlgo11,ApproxAlgo12}

\bibliographystyle{amsplain}
\bibliography{lotsize_optimization}

\label{lastpage}

\end{document}